
\documentstyle{amsppt}
\magnification 1200
\NoRunningHeads
\NoBlackBoxes
\document

\def\t{\hat t}
\def\Ua{U_q(\tilde\g)}
\def\U2{{\Ua}_2}
\def\a{\frak a}
\def\g{\frak g}
\def\h{\frak h}
\def\n{\hat n}

\def\Z{\Bbb Z}
\def\C{\Bbb C}

\def\Tr{\text{\rm Tr}}
\def\l{\lambda}

\def\<{\langle}
\def\>{\rangle}
\def\o{\otimes}
\def\e{\varepsilon}
\def\q{\hat q}

\def\Ug{U_q({\frak g})}

\topmatter
\title Central elements for quantum affine algebras and affine Macdonald's
operators
\endtitle
\author {\rm {\bf Pavel I. Etingof} \linebreak
\vskip .1in
Department of Mathematics\linebreak
Harvard University\linebreak
Cambridge, MA 02138, USA\linebreak
e-mail: etingof\@math.harvard.edu}
\endauthor

\endtopmatter
\centerline{December 19, 1994}
\vskip .1in
\centerline{\bf Abstract}
\vskip .1in
We describe a generalization of Drinfeld's description
of the center of a quantum group to the case of quantum affine algebras.
We use the obtained central elements to construct the affine analogue
of Macdonald's difference operators.
\vskip .1in

\centerline{\bf 1. The center of $U_q(\g)$, where $\g$ is a simple Lie
algebra.}

\vskip .1in
1.1. Let $\g$ be a simple Lie algebra over $\C$ of rank $r$.
Let $\h$ be a Cartan subalgebra in $\g$. Let $W$ be the Weyl group.
We fix a Weyl group invariant
inner product $\<,\>$ on $\h^*$ and $\h$ by setting
$\<\alpha,\alpha\>=2$ for short roots. Fix
a polarization of $\g$.
Let $\rho\in \h^*$ be the half-sum of positive roots.
Let $Q$ be the root lattice of $\g$, $Q^+$ be the semigroup
with $0$ spanned by the positive roots. Let $\alpha_1,...,\alpha_r$
be the simple roots. Let $A=(a_{ij})$,
$a_{ij}=2\<\alpha_i,\alpha_j\>/\<\alpha_i,\alpha_i\>$
be the Cartan matrix of $\g$.
Let $P$ be the weight lattice of $\g$, $P^+$ be the
set of dominant integral weights, and $\omega_1,...,\omega_r$ be the
fundamental weights. Let $N$ be the order of $P/Q$.
Note that for any $\l,\mu\in P$ $N\<\l,\mu\>\in\Z$.

1.2. Let $\q$ be a formal variable. Let $q=\q^N$. For any
$a\in\frac{1}{N}\Z$ we define $q^a:=\q^{Na}$.

Let $U_q(\g)$ be the
quantum group corresponding to $\g$
(\cite{Dr1,J1}). It is a Hopf algebra generated over the field $F=\C(\q)$ by
the elements $E_i$, $F_i$, $1\le i\le r$, and $K_\beta$, $\beta\in P$, freely
modulo the
relations:
$$
\gather
K_\beta K_\gamma=K_{\beta+\gamma},\ K_0=1,
\ K_\beta E_i K_{-\beta}=q^{\<\beta,\alpha_i\>}E_i,\\
\ K_\beta F_i K_{-\beta}=q^{-\<\beta,\alpha_i\>}F_i,
\ [E_i,F_j]=\delta_{ij}\frac{k_i-k_i^{-1}}{q_i-q_i^{-1}},\ q_i=
q^{\frac{1}{2}\<\alpha_i,\alpha_i\>},  k_i=K_{\alpha_i},\\
\sum_{l=0}^{1-a_{ij}}(-1)^l\biggl[\matrix 1-a_{ij}\\
l\endmatrix\biggr]_{q_i}E_i^{1-a_{ij}-l}
E_jE_i^{l-1}=0, i\ne j,\\
\sum_{l=0}^{1-a_{ij}}(-1)^l\biggl[\matrix 1-a_{ij}\\
l\endmatrix\biggr]_{q_i}F_i^{1-a_{ij}-l}
F_jF_i^{l-1}=0, i\ne j,
\endgather
$$
where $\biggl[\matrix m\\
l\endmatrix\biggr]_{q}:=\prod_{j=0}^{l-1}\frac{[m-j]_q}{[l-j]_q}$,
$[m]_q:=\frac{q^m-q^{-m}}{q-q^{-1}}$.

The comultiplication is defined by
$\Delta(E_i)=E_i\o k_i+1\o E_i$, $\Delta(F_i)=F_i\o 1+k_i^{-1}\o F_i$,
$\Delta(K_\beta)=K_{\beta}\o K_\beta$.
 The antipode is given by $S(E_i)=-E_ik_i^{-1}$, $S(F_i)=-k_iF_i$,
$S(K_\beta)=K_{-\beta}$. The counit $\e$ annihilates $E_i$ and $F_i$ and maps
$K_\beta$ to the
identity.

1.3.
Let $K$ be
an extension of $F$. We call any character $\theta: P\to K^*$ a $K$-weight
for $U_q(\g)$. We say that a representation $V$ of $U_q(\g)\o K$ is
with highest weight $\theta$ if it is generated by a vector $v$ such that
$E_iv=0$,
$1\le i\le r$, and $K_\beta v=\theta(\beta)v$, $\beta\in P$.
Note that any
 character $\theta: P\to K^*$ extends to a ring homomorphism $\C[P]\to K$. We
denote
this homomorphism also by $\theta$.

In particular, for $\mu\in P$ we can define
$\theta_\mu: P\to F^*$ by $\theta_\mu(\l)=q^{\<\mu,\l\>}$.

1.4.
Let $\Cal O_f$ be the category of finite-dimensional representations
$V$ of $U_q(\g)$ such that $V=\oplus_{\l\in P}V[\theta_\l]$,
$K_\beta|_{V[\theta_\l]}=
q^{\<\beta,\l\>}Id$. For $\mu\in P^+$, let
$V_\mu\in\Cal O_f$ be the irreducible representation with highest weight
$\theta_\mu$.

1.5. Let $U_q({\frak n^\pm})$ be the
subalgebras in $U_q(\g)$ generated by $\{E_i\}$,
$\{F_i\}$, respectively. Let $\{a_i,i\ge 0\}$ be a homogeneous basis
of $U_q({\frak n^+})$ ($a_0=1$), and let $\{a_i^*\}$ denote the dual basis to
$\{a_i\}$
with respect to the Drinfeld pairing between $U_q({\frak n^+})$ and
$U_q({\frak n^-})$ \cite{Dr1}.

1.6. Let $V\in \Cal O_f$. Define an element $C_V\in U_q(\g)$ by
$$
C_V:=\sum_{i,j,\l}
q^{2\<\l,\rho\>}\text{Tr}|_{V[\theta_\l]}(a_ia_j^*)a_i^*K_{2(\l-\beta_j)}a_j,
\tag 1
$$
where $\beta_i$ is the weight of $a_i$, i.e. $K_\gamma
a_iK_{-\gamma}=q^{\<\gamma,
\beta_i\>}a_i$, $\gamma\in P$.
Clearly, the sum in (1) is finite:
almost all terms vanish.

{\bf Remark. } The element $C_V$ can be shortly written in the form
$$\text{Tr}|_V(R_{21}R(1\o K_{2\rho})),$$ where
$R$ is the universal R-matrix for $U_q(\g)$ \cite{Dr1}, and $R_{21}$ is the
result
of permutation of the two components of $R$.

1.7.
\proclaim{Theorem 1}\cite{Dr1},\cite{R}

(i) $C_V$ belongs to the center $\Cal Z$ of $U_q(\g)$.

(ii) (quantum Harish-Chandra's theorem)
Let $\text{Rep}\g$ be the Grothendieck ring of the category $\Cal O_f$
as a tensor category. Then the map $V\to C_V$ is a ring isomorphism
$F\o_\Z\text{Rep}\g\to \Cal Z$.
Thus, $\{C_{V_\mu},\mu\in P^+\}$ is a basis of $\Cal Z$, and
$\Cal Z=F[C_{V_{\omega_1}},...,C_{V_{\omega_r}}]$.

(iii) (formula for the Harish-Chandra homomorphism) Let $K$ be an extension of
$F$ and $M$ be a highest weight representation of $U_q(\g)\o K$ with highest
weight $\sigma:P\to K^*$. Then $C_V|_M=\theta(\chi_V)$, where
$\theta=\sigma^2\theta_{2\rho}$,
and $\chi_V\in\C[P]$ is the character of $V$.

\endproclaim
\vskip .1in
\centerline{\bf 2. Generalization to the affine case.}
\vskip .1in

2.1. Let $\tilde \g=\g\o\C [t,t^{-1}]\oplus\C c\oplus\C d$ be the affine
Lie algebra corresponding to $\g$.
Let $\tilde\h=\h\oplus \C c\oplus\C d$ be the Cartan subalgebra in $\g$. Then
$\tilde\h^*=\h^*\oplus \C \e\oplus\C \delta$, so that
$\e(c)=\delta(d)=1$, and $\delta|_{\h\oplus\C c}=\e|_{\h\oplus\C d}=0$.
We can extend the form $\<,\>$ to a nondegenerate
form on $\tilde \h^*$ and $\tilde\h$ by setting
$\<\e,\delta\>=1$, and
$\<\e,\e\>=\<\delta,\delta\>=\<\e,\h^*\>=\<\delta,\h^*\>=0$.
Define $\alpha_0=\delta-\tilde\alpha$, where $\tilde\alpha$ is the maximal
root of $\g$. Then the simple roots of $\tilde\g$ are $\alpha_0,...,\alpha_r$.
The affine Weyl group $\hat W$ is generated by the finite Weyl group
$W$ and an additional reflection $s_0$ acting by
$s_0\l=\l-\<\l,\alpha_0\>\alpha_0$.
It acts in a natural way on $\tilde \h^*$.
This action preserves $\<,\>$.

Let the weight $\hat\rho\in \tilde\h^*$ be defined by the condition
$\<\hat\rho,\alpha_i\>=\frac{1}{2}\<\alpha_i,\alpha_i\>$
for all simple roots $\alpha_i$, $i=0,...,r$, and $\<\hat\rho,\e\>=0$.
Then $\hat\rho=\rho+g\e$, where $g$ is the dual Coxeter number of $\g$.

Let $\hat\omega_j\in\tilde \h^*$, $0\le j\le r$, be the fundamental
weights of $\tilde\g$, defined by $\<\hat\omega_j,\alpha_i\>=\frac{\delta_{ij}}
{2}\<\alpha_i,\alpha_i\>$, and $\<\hat\omega_j,\e\>=0$.
Define $\tilde P$ to be the set of all weights of the form
$\l=m\delta+\sum_{j=0}^rm_j\hat\omega_j$, $m,m_j\in\Z$.
Let $\tilde P^+$ be the set of dominant integral weights, i.e.
weights of the form
$\l=m\delta+\sum_{j=0}^rm_j\hat\omega_j$, $m\in\Z,m_j\in\Z_+$.
Let $\tilde Q$ be the root lattice of $\tilde\g$, and let $\tilde Q^+$ be the
semigroup with $0$ spanned by the positive roots.

2.2. Let $\Ua$ be the quantum
affine algebra corresponding to $\tilde\g$ (\cite{Dr1}). It is defined as
follows.
Its generators (over $F$) are $E_i,F_i$, $i=0,...,r$, $K_\beta$, $\beta\in
\tilde P$,
which satisfy exactly the same relations as in Section 1, with the only
difference
that $i,j$ vary from $0$ to $r$, and $\beta$ varies over $\tilde P$.
Let $U_q(\hat\g)$ be the subalgebra in $\Ua$ generated by $E_i,F_i,K_\beta$,
$\beta\in\tilde P$, $\<\beta,\delta\>=0$.

 Let $\{e^\mu,\mu\in \tilde P\}$ be the natural basis of $\C[\tilde P]$.
Let $F[[\tilde P]]_T$ denote the ring of formal series
$\sum_{m=0}^\infty c_me^{\l_m}$, $c_m\in F$,
$\l_m\in\tilde P$, such that
$\<\l_m,\hat\rho\>\to -\infty$, $m\to \infty$, and the order of the pole of
the function $c_m$ at $\q=0$ and $\q=\infty$ is bounded from above
as a function of $m$.
We call such formal series {\it tempered}
(this is why we use the subscript T in the notation).
Let $A$ be the intersection of all the translates
of $F[[\tilde P]]_T$ under the action of the affine Weyl group.
(It is easy to see that $\hat W$ does not preserve $F[[\tilde P]]_T$.)
For $\phi=\sum c_m
e^{\l_m}\in A$, we formally set $\xi(\phi)=\sum c_mq^{\<\l_m,\hat\rho\>}
K_{\l_m}$.

Let $U_q({\frak \n^\pm})$ be the
subalgebras in $U_q(\tilde\g)$ generated by $\{E_i\}$,
$\{F_i\}$, respectively. Let $\{a_i,i\ge 0\}$ be a homogeneous basis
of $U_q({\frak \n^+})$ ($a_0=1$), and let $\{a_i^*\}$ denote the dual basis to
$\{a_i\}$
with respect to the Drinfeld pairing between $U_q({\frak \n^+})$ and
$U_q({\frak \n^-})$.

Let $\widehat{\Ua}$ denote the completion of $\Ua$ consisting
of all formal series of the form
$\sum_{i,j}a_i^*\xi(\phi_{ij})a_j$, where $\phi_{ij}\in A$, and
the number of distinct weights of terms in the sum is finite (so
that $\widehat{\Ua}$ is a $\tilde P$-graded algebra). The algebra structure on
$\widehat{\Ua}$ is obvious.

2.3. Let $K$ be an extension of $F$. We call any character $\theta: \tilde P\to
K^*$ a
$K$-weight
for $U_q(\tilde\g)$. We say that a representation $V$ of $U_q(\g)\o K$ is
with highest weight $\theta$ if it is generated by a vector $v$ such that
$E_iv=0$,
$0\le i\le r$, and $K_\beta v=\theta(\beta)v$, $\beta\in \tilde P$.
In particular, for $\mu\in \tilde P$ we can define
$\theta_\mu: \tilde P\to F^*$ by $\theta_\mu(\l)=q^{\<\mu,\l\>}$.

Now assume that $K=L(s)$, where $s$ is a transcendental element over a field
$L\supset F$,
and $\theta(\alpha_j)\in L[s^{-1}]$,
$\theta(\delta)\in s^{-1}L[s^{-1}]$. We call such a weight admissible.
Then the map $\theta$ extends to an algebra homomorphism
$A\to L((s))$. Indeed, if $a=\sum c_me^{\l_m}\in A$ then
$\<\l_m,\e\>\to-\infty$, $m\to\infty$, so $\theta(a)$ converges in $L((s))$.
Note that the set of admissible weights is invariant under shifting by roots of
$\tilde\g$.

Moreover, if $\theta=\theta_\mu$, $\mu\in \tilde P$, and $\<\mu,\delta\>>0$,
then
$\theta$ extends to a ring homomorphism $A\to F_-$, $F_-=\C((\q^{-1}))$.
It follows from the fact that for such $\mu$ we have $\<\mu,\l_m\>\to -\infty$,
$m\to\infty$.
Similarly, if $\<\mu,\delta\><0$ then $\theta$ extends to a ring homomorphism
$A\to F_+$,
$F_+=\C((\q))$. For convenience we denote by $\tilde P(\pm)$ the subsets of
$\tilde P$ consisting of all $\mu\in\tilde P$
such that $\<\mu,\delta\>\in\Bbb N^{\pm}$.

2.4. Let $M$ be a module over $\Ua\o K$, where $K$ is an extension of $F$. For
any weight $\theta: \tilde P\to K^*$ we denote by $M[\theta]$ the space
of all vectors of weight $\theta$ in $M$. We say that $M$ has weight
decomposition if
$M=\oplus_\theta M_\theta$. We say that $M$ is smooth if it has a weight
decomposition,
and $U_q(\frak \n^+)v$ is a finite-dimensional space for any $v\in M$. Let
$\Cal S_K$ be the
category of smooth modules. Let $\Cal S_L^a$ be the
full subcategory in $\Cal S_{L(s)}$ of smooth modules having only admissible
weights,
and $\Cal S_{\pm}$ be the full subcategory of $\Cal S_F$ consisting
of modules having only
weights of the form $\theta_\mu$, $\mu+\hat\rho\in
\tilde P(\pm)$. Then, as follows from
Section 2.3, for any $M\in S_L^a$,
$M\o L((s))$ naturally extends to a module over $\widehat{\Ua}\o L((s))$, and
for any
$M\in \Cal S_\pm$, $M\o F_{\mp}$ naturally extends to a module over
$\widehat{\Ua}\o F_{\mp}$. The action of $\xi(A)$ in these extensions are
defined by
$\xi(\phi)|_{M[\theta]}=\theta(\phi)\theta_{\hat\rho}(\phi)\text{Id}$.

2.5. We would like to construct some central elements in $\widehat{\Ua}$.
One of these elements is obvious -- it is $K_{\delta}$. The rest of the
central elements can be constructed analogously
to Drinfeld's method.

Let $\U2$ be the subalgebra in $\Ua$ generated by $U_q(\hat\g)$ and $K_\e^2$.
It is clear that $\U2$ is a Hopf subalgebra, and $\Ua=\U2\oplus K_\e\U2$.
Let $\tilde P_2$ be the sublattice of $\tilde P$ spanned by $\hat\omega_i$ and
$2\delta$,
and let $\tilde P_{1/2}$ be the
lattice of $\tilde P$ spanned by $\hat\omega_i$ and $\frac{1}{2}\delta$.
If $\l\in P_{1/2}$, we define the character $\theta_\l:P_2\to F^*$ by
the usual formula $\theta_\l(\beta)=q^{\<\l,\beta\>}$.

Let $\Cal O_i$ be the category of integrable representations
$V$ of $\U2$ such that $V=\oplus_{\l\in \tilde P_{1/2}}V[\theta_\l]$,
$\text{dim}V[\theta_\l]
<\infty$, and
$K_\beta|_{V[\theta_\l]}=
q^{\<\beta,\l\>}Id$, $\beta\in\tilde P_2$. (Integrability means that the
restriction of $V$
to any $U_q(sl_2)$-subalgebra corresponding to a simple root is a direct sum of
finite dimensional modules).
This category is semisimple (\cite{L}): any object is a (possibly infinite)
direct sum of
irreducible objects.

For $\mu\in P_{1/2}$, let
$V_\mu\in\Cal O_i$ be the irreducible representation of $\U2$
with highest weight
$\theta_\mu$. Any simple object in $\Cal O_i$ is isomorphic to $V_\mu$ for a
unique $\mu$.
Characters of representations $V\in \Cal O_i$ are infinite linear combinations
of $e^\mu$, $\mu\in \tilde P_{1/2}$.

Clearly, $\Cal O_i$ is closed with respect to the tensor product.

2.6.
Let $V\in \Cal O_i$. Define an element $C_V\in U_q(\tilde\g)$ by
$$
C_V:=\sum_{i,j,\l}
q^{2\<\l,\hat\rho\>}\text{Tr}|_{V[\theta_\l]}(a_ia_j^*)a_i^*K_{2(\l-\beta_j)}a_j,
\tag 2
$$
where $\beta_i$ is the weight of $a_i$. Clearly, the sum in (2) is a well
defined
element of $\widehat{\Ua}$. For example, if $V_n$ is the 1-dimensional
$\U2$-module in which $U_q(\hat\g)$ acts trivially and
$K_\e^2$ acts by multiplication by $q^n$ then $C_{V_n}=q^{gn}K_{\delta}^{n}$.

{\bf Remark. } 1. The element $C_V$ can be shortly written in the form
$$\text{Tr}|_V(R_{21}R(1\o K_{2\hat\rho})),$$ where
$R$ is the universal R-matrix for $U_q(\tilde\g)$ \cite{Dr1}, and $R_{21}$ is
the result
of permutation of the two components of $R$.

2. The reason we need to introduce the algebra $\U2$ is the following.
If we only consider elements $C_V$ with $V$ being an integrable $\Ua$-module,
like we did in the finite-dimensional case, then the algebra generated
by all such element will not contain $K_\delta$
(although it will contain its square).

2.7. The following theorem is the affine analogue of Drinfeld's theorem.

\proclaim{Theorem 2} (i) $C_V$ belongs to the center $\tilde\Cal Z$ of
$\widehat{\Ua}$.

(ii) Let $\text{Rep}\tilde\g$ be the Grothendieck ring of the category
$\Cal O_i$
as a tensor category. Then the map $V\to C_V$ is a ring isomorphism
$F\o_\Z\text{Rep}\tilde\g\to \tilde\Cal Z$.
Thus, $\{C_{V_\mu},\mu\in \tilde P^+\}$ is a topological basis of
$\tilde\Cal Z$, and
$\tilde\Cal Z=F[C_{V_{\hat\omega_0}},...,C_{V_{\hat\omega_r}}]((K_{\delta}))_T$
(here, as before, the subscript T means the tempered series, i.e. such that the
order of poles of
of its coefficients at $\q=0$ and $\q=\infty$ are bounded from above).

(iii) Let $L$ be an extension of $F$, $K=
L(s)$ be the field of rational functions of $s$.
Let $M$ be a highest weight representation of $U_q(\tilde\g)\o K$ with an
admissible highest weight $\sigma: \tilde P\to K^*$.
Then $M\o L((s))$ extends to a representation of
$\widehat{\Ua}\o L((s))$, and $C_V$ acts there as $\theta(\chi_V)Id$, where
$\chi_V$ is the character of $V$, and $\theta=\sigma^2\theta_{2\hat\rho}$.

(iv) Let $\nu\in \tilde P$ and $\<\nu,\delta\>+g\in \mp\Bbb N$. Let $M$ be any
representation
of $\Ua$ with highest weight $\nu$. Then $M\o F_{\pm}$ extends to a
representation of $\widehat{\Ua}\o F_{\pm}$ in a natural way, and
$C_V$ acts there as $\theta_{2(\nu+\hat\rho)}(\chi_V)Id$, where
$\chi_V$ is the character of $V$.

\endproclaim

2.8. {\bf Remarks. } 1. The statement about topological basis in (ii) means
that
any element $z\in\tilde\Cal Z$ can be uniquely represented as
$\sum_nb_nC_{V_{\mu_n}}$,
$b_n\in F$, and $\<\mu_n,\hat\rho\>\to -\infty$ as $n\to\infty$.

2. We formulate (iv) separately from (iii) because weights in $\tilde P$ are
not
admissible.

3. This theorem can be easily generalized to the case of $U_q(\a)$ --
quantization of an arbitrary symmetrizable Kac-Moody algebra $\frak a$. The
changes which
need to be made are obvious, and we do not discuss them here.

4. In the classical case $q=1$, Kac \cite{K} constructed central elements in a
completion
of $U(\a)$, where $\a$ an arbitrary symmetrizable Kac-Moody algebra. Kac's
construction assigns such a central element to any function on the Tits cone
in the set of weights, in such a way that the this assignment inverts the
Harish-Chandra
homomorphism. Kac's construction easily generalizes to quantum groups, and the
elements
$C_V$ are examples of central elements one can obtain using this construction.
But in the quantum case we also have nice explicit formulas like (1),(2), which
are
absent in the classical case.
\vskip .1in

\demo{2.9. Proof of Theorem 2} The proof is analogous to the proof of Theorem
1.

(i) Let $V_1,V_2$ be highest weight representations of
$\Ua$, and let $X:V_1\o V_2\to V_1\o V_2$ be an intertwining operator. Assume
that the trace $\text{Tr}|_{V_2}(X(1\o K_{2\hat\rho}))$ makes sense. Then
it defines an intertwining operator $V_1\to V_1$, i.e. commutes with the action
of $\Ua$. Therefore, (i) follows from Remark 2.6
and the fact that $R_{21}R$ commutes with $\Delta(a)$,
$a\in\Ua$.

(iii),(iv).
When we apply
$C_V$ to the highest weight vector
$v\in M$, all terms in (2) with $i,j>0$ vanish, and what remains is
$$
C_Vv=\sum_\l
q^{\<2\l,\hat\rho\>}\text{dim}V[\theta_\l]K_{2\l}v=\theta(\chi_V)v.
$$
Since $M$ is generated by $v$, we get (iii),(iv).

(ii) The fact that the assignment $V\to C_V$ is a ring homomorphism follows
from
the fusion property of the $R$-matrix. Moreover, Kac \cite{K} showed that a
central element $C\in\widehat{\Ua}$ of the form $\sum a_i^*\xi(\phi_{ij})a_j$
is completely determined by $\phi_{00}$ (strictly speaking, Kac's paper
treats the case $q=1$, but the results we need are obviously
valid in the quantum case).
Thus, the
map $V\to C_V$ is injective. To show that it is surjective, let us recall
\cite{K}
that $\theta_\l(\phi_{00})=\theta_{s_\alpha(\l)}(\phi_{00})$
when $\<\l+\hat\rho,\alpha\>=\frac{n}{2}\<\alpha,\alpha\>$, where
$n\in\Bbb N$, $\alpha$ is a positive
root of $\tilde \g$, and $s_\alpha$ is the Weyl reflection
with respect to
$\alpha$ (the Kac-Kazhdan condition); this is because the Verma module with
highest weight
$\theta_\l$
contains the Verma module with highest weight
$\theta_{s_\alpha(\l+\hat\rho)-\hat\rho}$
when the Kac-Kazhdan condition is satisfied, and $C_V$ acts by
$\theta_{\l+\hat\rho}(\phi_{00})Id$
on the Verma module with highest weight $\theta_\l$.
Therefore, $\phi_{00}\in A^{\hat W}$, where $A^{\hat W}$ is the set of $\hat
W$-invariant elements in $A$.
But the algebra $A^{\hat W}$ is nothing else but the completed algebra of
invariant
theta functions associated to the affine root system of $\tilde\g$.
By a theorem of Looijenga (\cite{Lo}), this algebra
is topologically spanned over $F$ by the characters of integrable
representations, so we get (ii). $\square$

{\bf Remarks.} 1. As we have seen, the construction of the center for
$U_q(\hat\g)$ given above works at any non-critical value of the central
charge: $K_{\delta}\ne q^{-g}$. At the critical value $q^{-g}$, the structure
of
the center is different. In fact, in this case the center is much
bigger -- it contains an algebra of polynomials of infinitely many variables.
However, even in this case the above construction (formula (2))
allows to obtain a central element in a completion of $U_q(\hat\g)$
if the module $V$ is taken not from the category $\Cal O_i$ but from
the category of finite-dimensional $U_q(\hat\g)$-modules.
This construction was discovered in \cite{RS}; see also \cite{DE}.

2. The constructions of central elements
for $U_q(\g)$, $U_q(\hat\g)$ by formulas (1),(2), as well as the construction
of central elements of $U_q(\hat\g)$ at the critical level
given in \cite{RS},\cite{DE}, mentioned in the previous remark,
are special cases of the general construction of categorical trace
defined by J.Bernstein \cite{B}, as follows.

Let $A,B$ be categories, $F:A\to B$ be a functor, and $F^*,^*F:B\to A$
be right and left adjoint functors to $F$. Let $\e: F^*\to ^*F$ be a
morphism of functors. Let $X$ be an object in $A$ and $a:FX\to FX$
be an endomorphism. Consider the morphism $tr(a):X\to X$ defined as
the composition of morphisms
$tr(a)=i_X\circ ^*F(a)\circ \e_{FX}\circ j_X$, where
$i_X: F^*FX\to X$, $j_X: X\to ^*FFX$ are the adjunction morphisms.
This defines the linear operator of categorical trace
$tr: \text{End}(F)\to \text{End}(Id_A)$, where $Id_A:A\to A$ is the
identity functor.

Let $U$ be an associative algebra over a ring $\Cal F$, and $A$ be a full
subcategory of the category of $U$-modules.
Then the center $Z(U)$ of $U$ naturally maps to the ring $\text{End}(F)$,
where $F$ is any functor on $A$. In particular, if $F$ is the identity
functor, this map is often an embedding.
Therefore $Z(U)$ can often be identified with
a subring of the ring of endomorphisms of the identity functor on the category
of representations of $U$.

It is shown in \cite{B} how to use the trace construction to produce many
central elements of $U$ when we have only
 one fixed central element $C\in Z(U)$.
Indeed, let $F: A\to A$ be any functor satisfying the above conditions.
Then we can consider $tr F(C)\in \text{End}(Id_A)$, and if we are lucky,
so that this element is in the image of the center, then we get a new central
element $C_F\in Z(U)$. Let us demonstrate how this works
for quantum groups.

If $U$ is a qausitriangular Hopf algebra
 (e.g. $U_q(\g), U_q(\hat\g)$) then one can define
an element $u=m((S\o 1)(R_{21}))$
in a completion of this algebra,
where $m$ is the multiplication map,
$S$ is the antipode, and $R_{21}$ is the universal $R$-matrix with permuted
components. Conjugation by $u$ is $S^2$.
This element was defined by Drinfeld \cite{Dr2}.
It is easy to show that the element $Z=uS(u)$ is central. Moreover,
in the cases of $U_q(\g)$, $U_q(\hat\g)$ it is possible to find
a central element $C$ such that $\Delta(C)=(C\o C)(R_{21}R)$,
and $C^2=Z^{-1}$
(it is obtained by formal extraction of square root from $Z^{-1}$;
Hopf algebras where $C$ exists are called ribbon Hopf algebras, \cite{RT}).
Now, if $V$ is an irreducible,
integrable highest weight module, and $F=F_V$ is the functor
of tensoring of $U$-modules with $V$ (on the right). Then
$^*F, F^*$ are the functors of tensoring with $^*V,V^*$, respectively,
and $u: V^*\to ^*V$ defines an isomorphism $\e: F^*\to ^*F$.
It is easy to see that $C_V=\frac{1}{C|_V}C^{-1}C_F$, where $C_F$ is the
categorical trace of $C$, and $C|_V$ is the eigenvalue of $C$ in $V$.
This explains formulas (1),(2).

\vskip .1in
\centerline{\bf 3. Central elements and Macdonald operators}
\vskip .1in

3.1. Macdonald operators acting on the space of rational functions of
$n$ variables over $\C(q,t)$ are defined as follows:
$$
M^r=t^{r(r-n)} \sum\limits_{i_1<i_2<\ldots<i_r}
  \biggl(\prod\Sb j\notin\{i_1\ldots i_r\}\\ l=1\ldots r \endSb
       \frac{t^2x_{i_l} -x_j}{x_{i_l}-x_j} \biggr)
	T_{i_1}\ldots T_{i_r}\tag 3$$
where $(T_{i}f)(x_1,\ldots,x_n)=f(x_1,\ldots, q^2x_i,\ldots,
x_n)$, $t$ is a formal variable, and $r=1,\ldots,n-1$.

Macdonald proved \cite{M} that these operators are pairwise
commutative: $[M^r,M^s]=0$, $1\le r,s\le n-1$.

In \cite{EK1},\cite{EK2} it was shown how to obtain Macdonald
operators from central elements of the quantum group $\Ug$, where
$\g=\frak{sl}_n$. This was done as follows.

3.2. Recall that fundamental representations of $U_q(\g)$ are
$q$-deformed exterior powers $(\Lambda^rF^n)^*$, $r=1,...,n-1$.
Consider the elements
$$
C_r:=C_{(\Lambda^rF^n)^*}\in\Ug,\tag 4
$$
where $\rho=(\frac{n-1}{2},\frac{n-3}{2},...,\frac{-n+1}{2})$.
These elements are central and freely generate the center of $\Ug$,
according to Theorem 1.

3.3. Let $\t$ be a formal variable such that $\t^2=t$.
Let $K$ be an extension of $F(\t)$, and let $\theta: P\to K^*$ be a weight.
Let $M_{\theta}$ denote the Verma module over $\Ug$ with highest weight
$\theta$, and let
$$
U=\{p\in K[x_1^{\pm 1},\ldots,
x_n^{\pm 1}], \text{deg }p=0\}\tag 5
$$
with the action of $\Ug\o K$ given by
$$\gathered
E_i\mapsto x_i D_{i+1},\
F_i\mapsto x_{i+1} D_i, \\ K_\beta f(x_1,...,x_n)=
f(q^{\beta_1}x_1,...,q^{\beta_n}x_n),\beta=(\beta_1,...,\beta_n)\in P\\
(D_i f)(x_1,\ldots, x_n)=\frac{tq^{-1}f(x_1,\ldots, qx_i,\ldots ,x_n)-
t^{-1}qf(x_1,\ldots, q^{-1}x_i,\ldots ,x_n)}{(q-q^{-1})x_i}
\endgathered\tag 6$$

The set of
weights of $U$ is the root lattice $Q$, and every
weight subspace is one-dimensional:
$U[\mu]=K x_1^{\mu_1}...x_n^{\mu_n}$, $\mu=(\mu_1,...,\mu_n)\in Q$.
We fix a vector $u\in U[0]$.

3.4. \proclaim{Lemma 1} The following two conditions on the Verma
module $M_\theta$ are equivalent:

(i) $M_\theta$ is irreducible.

(ii) There exists a
unique intertwiner
$$\Phi_\theta\colon M_{\theta}\to M_{\theta}
\hat\otimes U,\tag 7 $$
where $\hat\otimes$ denotes the completed tensor product with respect to the
principal
grading in $M_\theta$, normalized by the condition
$\Phi_\theta v=v\o u+\text{lower order terms}$, where $v$ is the highest weight
vector of $M_\theta$.
\endproclaim

Proof of the (i)$\to$ (ii) part is based on the general fact:
if a Verma module $M_\theta$  is irreducible
and $U$ is any $U_q(\g)\o K$-module
then the space of intertwiners $M_\theta\to M_\theta\otimes U$ is in one-to-one
correspondence with the zero-weight subspace $U[0]$, i.e the space of
$K_\beta$-invariant vectors
in $U$ (for all $\beta\in P$).
The (ii)$\to$ (i) part (which we don't use in this paper)
follows from the study of singularities of
$\Phi_\theta$ as a function of $\theta$ done in \cite{ES}. (In \cite{ES},
the case $q=1$ is treated, but in the quantum case the situation
is analogous).

{\bf Remark.} The module $M_\theta$
over $U_q(\g)$ is reducible if and only if the quantum Kac-Kazhdan condition
is satisfied: $\theta(2\alpha)=q^{n\<\alpha,\alpha\>}$ for some positive
root $\alpha$ of $\g$, and some integer $n>0$. This follows from the formula
for the determinant of the Shapovalov form, see \cite{DCK}.

3.5. Let $Y=\text{Hom}(P,K^*)$. Define a completion $K[[Y]]$ of the group
algebra $K[Y]$ of $Y$.
By definition, $K[[Y]]$ is spanned by the subspaces $K[[Y]]_\sigma, \sigma\in
Y$, where
$K[[Y]]_\sigma$ consists of all, possibly infinite, formal series of the form
$\sum_{\beta\in Q^+}c_\beta\sigma\theta_{-\beta}$, where $c_\beta\in K$.

Let $M$ be a highest weight $U_q(\g)\o K$-module.
For any $\Phi: M\to M\hat\o U$ denote by $\Phi[\theta]$ the diagonal block of
$\Phi$ of weight $\theta$: $\Phi[\theta]: M[\theta]\to M[\theta]\o U[0]$.
Let $K(Y)$ be the quotient field of $K[Y]$, and let
$X:M\o K(Y)\to M\o K(Y)$ be the linear operator defined by
$Xv=\theta v$, $v\in M[\theta]$.

Let $\tau_0:P\to K^*$ be defined by $\tau_0(\l)=\t^{2\<\l,\rho\>}$, and let
$\tau=\tau_0\theta_{-\rho}$.

Set
$$
\Psi_{\theta}=\Tr|_{M_{\theta\tau}\o K(Y)}(\Phi_{\theta\tau}X):=
\sum_{\beta\in Q^+}\Tr(\Phi_{\theta\tau}[\theta\tau\theta_{-\beta}])
\theta\tau\theta_{-\beta},\tag 8
$$
This is an element of $K[[Y]]\o U[0]$.
The construction of $\Psi_\theta$ is valid for a generic
$\theta\in Y$, i.e. such that the corresponding Verma module $M_{\theta\tau}$
is
irreducible.

Since the space $U[0]$ is one-dimensional,
we can
identify it with $K$ by sending $u$ to $1$ and
regard $\Psi_\theta$ as a series with scalar coefficients, i.e. as an element
of
$K[[Y]]$.

3.6. Now we explain a connection between Macdonald difference operators
and the series $\Psi_\theta$.

\proclaim{Definition} A standard difference operator is any
Laurent polynomial in the operators $T_{i}$ whose
coefficients are rational functions of $y_1=x_2/x_1,...,y_{n-1}=x_n/x_{n-1}$
(over $F(\hat t)$)
which are regular at $y_1=...=y_{n-1}=0$.
\endproclaim

It is clear that the algebra of standard difference operators can be embedded
in
the algebra $K[[y_1,...,y_{n-1}]][T_1,...,T_n]$
(by taking the Taylor expansion of the coefficients).
It therefore
follows that any standard difference operator
naturally acts on the space $K[[Y]]_\sigma$, $\sigma\in Y$.
This action is given by $T_{j}\theta=\theta(\omega_j-\omega_{j-1})^2\theta$
(by convention $\omega_0=0$), and
$y_j\theta=\theta\theta_{-\alpha_j}$.
There is a natural power series topology on $K[[Y]]_\sigma$, and the action of
any
standard difference operator is continuous in this topology.

3.7. \proclaim{Theorem 3}(\cite{EK2})

(i) For every $r$ between $1$ and $n-1$
there exists a unique standard difference operator $\tilde M^r$
defined over $F(t)$, such that
$$
\Tr|_{M_{\theta\tau}\o K(Y)}(\Phi_{\theta\tau}C_rX)=
\tilde M^r\Tr|_{M_{\theta\tau}\o K(Y)}(\Phi_{\theta\tau}X)\tag 9
$$
for all $\theta$ not satisfying the Kac-Kazhdan condition.

(ii) $[\tilde M^r,\tilde M^s]=0$, $1\le r,s\le n-1$.

(iii) Introduce the series
$$
\varphi=\tau\prod_{i=1}^{\infty}\prod_{\alpha\in R^+}
\frac{1-q^{2i}\theta_{-\alpha}}{1-q^{2(i-1)}t^2\theta_{-\alpha}}
\in \text{Aut}K[[Y]],\tag 10
$$
where $R^+$ is the set of positive roots of $\g$, and
$\tau$ denotes the operator of multiplication by $\tau$ in $K[[Y]]$.
Then $\frac{1}{\varphi}\tilde M^r(\varphi f)=M^rf$ for any
$f\in K[[Y]]$, where $M^r$ are the Macdonald operators defined
by (3). Thus, for any $\theta\in Y$ for which $\Psi_{\theta\tau}$ is defined
the series $\psi_{\theta}:=\Psi_{\theta}/\varphi$,
is a common eigenvector
of the Macdonald operators in $K[[Y]]$:
$M^r\psi_\theta=P_r(\theta)\psi_\theta$,
where $P_r(\theta)=(\theta\tau_0)^2(\chi_r)$, and $\chi_r$
is the character of $(\Lambda^rF^n)^*$.

(iv) Let $\sigma\in Y$ be such that the operator $\Phi_{\theta\tau}$
exists for any
$\theta=\sigma\theta_{-\beta}$, $\beta\in Q^+$. Then the vectors
$\psi_{\sigma\theta_{-\beta}}$, $\beta\in Q^+$, form a topological basis of
 $K[[Y]]_\sigma$  (with respect to the usual power series topology).
Thus, $M^r$ are simultaneously diagonalizable
in $K[[Y]]_\sigma$.  Generically, each of them has a simple spectrum
in $K[[Y]]_\sigma$.

(v) Let $\sigma\in Y$ be a weight such that the weights
$\sigma^w:=w(\sigma\tau_0)\tau_0^{-1}$, $w\in S_n$, are
all distinct, and $\Phi_{\sigma^w\tau}$ exists for each $w$.
 Then the system of difference equations $M^r\psi=P_r(\sigma)\psi$,
$r\in 1,...,n-1$, has $n!$ linearly independent solutions in $\Cal K[[Y]]$.
They are $\psi_{\sigma^w}$, $w\in S_n$.

\endproclaim

3.8.
In particular, the construction of functions $\psi_\theta$ can be specialized
to the
case when $\theta=\theta_\l$, $\l$ is a weight in $P^+$. Lemma 1 implies that
the operator
$\Phi_{\theta\tau}$ always exists in this case, and
we have

\proclaim{Theorem 4} If $\l\in P^+$ then the function $\psi_{\theta_\l}$
belongs
to $\C(q,t)[P]$ and coincides with the Macdonald polynomial
$P_\l(x;q,t)$ (see \cite{EK2}).
\endproclaim

\vskip .1in

\centerline{\bf 4. Affine Macdonald operators.}
\vskip .1in

In this section we will define an affine analogue of Macdonald operators using
the analogy with Section 3. We consider the case $\g=sl(n)$.

4.1. Let $p:U_q(\hat\g)\to \Ug$ be the Jimbo homomorphism
\cite{J2}. Let $U(1)=p^*(U)$.

Let $\theta:\tilde P\to K^*$ be a weight, and let
$M_\theta$ be the Verma module over $\Ua\o K$
with highest weight $\theta$.

4.2. \proclaim{Lemma 2} The following two conditions on the Verma module
$M_\theta$ are equivalent.

(i) $M_\theta$ is irreducible.

(ii) There exists a
unique $U_q(\hat\g)$-intertwiner
$$\Phi_\theta\colon M_{\theta}\to M_{\theta}
\hat\otimes U(1), $$
where $\hat\otimes$ denotes the
completed tensor product with respect to the principal
grading in $M_\theta$, normalized by the condition
$\Phi_\theta v=v\o u+\text{lower order terms}$, where $v$ is the highest weight
vector of $M_\theta$.
\endproclaim

Proof is analogous to the finite-dimensional case.

{\bf Remark.} The module $M_\theta$
over $U_q(\hat\g)$ is reducible if and only if the Kac-Kazhdan condition
is satisfied: $\theta(2\alpha)=q^{n\<\alpha,\alpha\>}$ for some
root $\alpha>0$ of
$\tilde\g$ and some integer $n>0$.

4.3.
Let $\tilde Y=\text{Hom}(\tilde P,K^*)$.
Define a completion $K[[\tilde Y]]$ of the group algebra $K[\tilde Y]$ of
$\tilde Y$.
By definition, $K[[\tilde Y]]$
is spanned by the subspaces $K[[\tilde Y]]_\sigma, \sigma\in \tilde Y$, where
$K[[\tilde Y]]_\sigma$ consists of all, possibly infinite, formal series of the
form
$\sum_{\beta\in Q^+}c_\beta\sigma\theta_{-\beta}$, where $c_\beta\in K$.

Let $M$ be a highest weight $U_q(\hat\g)\o K$-module.
For any $\Phi: M\to M\hat\o U$ denote by $\Phi[\theta]$ the diagonal block of
$\Phi$ of weight $\theta$: $\Phi[\theta]: M[\theta]\to M[\theta]\o U[0]$.
Let $K(\tilde Y)$ be the quotient field of $K[\tilde Y]$, and let
$X:M\o K(\tilde Y)\to M\o K(\tilde Y)$ be the linear operator defined by
$Xv=\theta v$, $v\in M[\theta]$.

Let $\tau_0:\tilde P\to K^*$
be defined by $\tau_0(\l)=\t^{2\<\l,\hat\rho\>}$, and let
$\tau=\tau_0\theta_{-\hat\rho}$.

Set
$$
\Psi_{\theta}=\Tr|_{M_{\theta\tau}\o K(\tilde Y)}(\Phi_{\theta\tau}X):=
\sum_{\beta\in \tilde Q^+}\Tr(\Phi_{\theta\tau}[\theta\tau\theta_{-\beta}])
\theta\tau\theta_{-\beta},\tag 11
$$
This is an element of $K[[\tilde Y]]\o U[0]$.
The construction of $\Psi_\theta$ is valid for a generic
$\theta\in \tilde Y$,
i.e. such that the corresponding Verma module $M_{\theta\tau}$ is
irreducible.

Since the space $U[0]$ is one-dimensional,
we can
identify it with $K$ by sending $u$ to $1$ and
regard $\Psi_\theta$ as a series with scalar coefficients, i.e. as an element
of
$K[[\tilde Y]]$.

4.4.
\proclaim{Definition} An affine difference operator is a formal expression
of the form $M=\sum a_mT_{\nu_m}$, $\nu_m\in\tilde P$,
such that for any
$w\in\hat W$ $\<w\nu_m, \hat\rho\>\to -\infty$, $m\to\infty$, and
$a_m\in F(t)[[\tilde Y]]_1$, where $1$ denotes the identity character
$\tilde P\to \C^*$.
\endproclaim

Affine difference operators form an algebra, multiplication in which is
defined by $T_\nu\theta=\theta(\nu)^2\theta T_\nu$.
This formula also gives rize to a natural action of the algebra of
affine difference operators in $K[[\tilde Y]]_\sigma\o L((s))$ for any
admissible $\sigma\in \tilde Y$. Also, if the series
$M$ is tempered, in the sense of Section 2, the action of
$M$ is defined in $G_{\pm}[[\tilde Y]]_\sigma$, $\sigma\in \tilde P(\mp)$,
$G_{\pm}=\C(\t)((\q))$. These actions are continuous in the standard
power series topology. We will further assume that all highest weights we
consider
are either admissible or from $\tilde P(\pm)$.

4.5. \proclaim{Theorem 5}

(i) For every $r$ between $0$ and $n-1$
there exists a unique affine difference operator $\tilde M^{\hat\omega_r}$
such that
$$
\Tr(\Phi_{\theta\tau}C_{V_{\hat\omega_r}}X)=
\tilde M^{\hat\omega_r}\Tr(\Phi_{\theta\tau}X)\tag 12
$$
for all $\theta$ not satisfying the Kac-Kazhdan condition.

(ii) $[\tilde M^{\hat\omega_r},\tilde M^{\hat\omega_s}]=0$, $0\le r,s\le n-1$.
\endproclaim

The proof of this theorem is analogous to the proof of Theorem 3.
Uniqueness follows from the fact that an affine difference operator
is uniquely determined by its action in $K[[\tilde Y]]_\sigma$
for generic $\sigma$, and this action is defined uniquely by (12).
Indeed, suppose we have an expression $\sum a_mT_{\nu_m}$ which acts
by $0$ on $K[[\tilde Y]]_\sigma$ for a generic admissible
weight $\sigma$. In particular, applying it to
$\sigma$, we get $\sum a_m\sigma(\nu_m)^2=0$.
This identity for generic admissible
$\sigma$ implies immediately that $a_m=0$ for all $m$.

4.6. Introduce the series
$$
\varphi=\tau\prod_{i=1}^{\infty}\prod_{\alpha\in \tilde R^+}
\frac{1-q^{2i}\theta_{-\alpha}}{1-q^{2(i-1)}t^2\theta_{-\alpha}}\in
F(t)[[\tilde Y]]_1,\tag 13
$$
where $\tilde R^+$ is the set of positive roots of $\tilde\g$.

Introduce the affine difference
operators $M^{\hat\omega_r}$ defined
by $M^{\hat\omega_r}f=
\frac{1}{\varphi}\tilde M^{\hat\omega_r}(\varphi f)$ for any
$f\in K[[\tilde Y]]$. These operators are pairwise commutative.

\proclaim{Definition} We call the operators $M^{\hat\omega_r}$
affine Macdonald operators.
\endproclaim

4.7. For any $\theta$ for which $\Psi_{\theta\tau}$ is defined consider
the series $\psi_{\theta}:=\Psi_{\theta}/\varphi$.

\proclaim{Corollary} (i) $\psi_\theta$ is a common eigenvector
of the Macdonald operators in $K[[\tilde Y]]_\theta$:
$M^{\hat\omega_r}\psi_\theta=\tilde P_r(\theta)\psi_\theta$,
where $\tilde
P_r(\theta)=(\theta\tau_0)^2(\chi_{\hat\omega_r})$, and $\chi_{\hat\omega_r}$
is the character of the $r-th$ fundamental representation.

(ii) Let $\sigma\in \tilde Y$ be such that the operator $\Phi_{\theta\tau}$
exists for any
$\theta=\sigma\theta_{-\beta}$, $\beta\in \tilde Q^+$. Then the vectors
$\psi_{\sigma\theta_{-\beta}}$, $\beta\in \tilde Q^+$, form a topological basis
of
 $K[[\tilde Y]]_\sigma$  (with respect to the usual power series topology).
Thus, $M^{\hat\omega_r}$ are simultaneously diagonalizable
in $K[[\tilde Y]]_\sigma$. Generically, each of them has a simple spectrum
in $K[[\tilde Y]]_\sigma$.

(iii) Let $\sigma\in \tilde Y$ be a weight such that the weights
$\sigma^w:=w(\sigma\tau_0)\tau_0^{-1}$, $w\in \hat W$, are
all distinct, and $\Phi_{\sigma^w\tau}$ exists for each $w$.
 Then the space of solutions of
the system of difference equations $M^{\hat\omega_r}\psi=\tilde
P_r(\sigma)\psi$,
is topologically spanned by the linearly independent solutions
$\psi_{\sigma^w}$, $w\in \hat W$.

\endproclaim

The proof is analogous to the finite-dimensional case.

4.8.
In particular, the construction of functions $\psi_\theta$ can be specialized
to the
case when $\theta=\theta_\l$, $\l$ is a weight in $P^+$. Then
the weight $\theta\tau$ is admissible ($K=F(s), s=\hat t^{-1}$), and Lemma 2
 implies that the operator
$\Phi_{\theta\tau}$ exists. In this case, we have

\proclaim{Theorem 6}\cite{EK3}
If $\l\in P^+$ then the function $\psi_{\theta_\l}$ belongs
to $\C(\q,\t)\o_F A^{\hat W}$ and coincides with the affine Macdonald
polynomial
$\hat J_\l$ defined in \cite{EK3}.
\endproclaim
(In fact, it is obvious
that the coefficients of $\hat J_\l$ are in $\C(q,t))$.
As shown in \cite{EK3}, affine Macdonald polynomials are a basis of
$\C(\q,\t)\o_F A^{\hat W}$ triangular with respect to the basis of characters
$\chi_\l$. Therefore, we find:

\proclaim{Corollary} The operators $M^{\hat\omega_r}$ are $\hat W$-invariant.
\endproclaim

This follows from the fact that an affine difference operator is uniquely
defined
by its action on symmetric theta-functions.

Thus , analogously to the finite-dimensional case,
affine Macdonald operators
 act on the space of symmetric theta-functions and are triangular
with respect to the standard ordering of dominant weights.
Affine Macdonald polynomials can be defined as their common eigenbasis.

In particular, it follows from this corollary that the coefficients to all
$T_\nu$ in
the affine Macdonald operators are of the form $\sum c_m \theta_{-\alpha_m}$,
where $\<w\alpha_m,\hat\rho\>\to +\infty$, $m\to\infty$, for all $w\in\hat W$.
However, these coefficients need not be Weyl group invariant.

\vskip .1in

\centerline{\bf Acknowledgements}

I would like to thank H.Garland who explained to me the idea
of a (conjectured) different constuction of affine
Macdonald operators. I am grateful to P.Polo for
explaining to me the connection of paper \cite{B} with
construction of center of quantum groups. I am grateful to
J.Bernstein, I.Frenkel, A.Kirillov Jr., and D.Kazhdan for discussions.
This work was partially supported by an NSF postdoctoral fellowship.

\end